\documentclass[12pt,leqno]{amsart}
\usepackage{float}
\usepackage{amssymb,amsmath,amsthm,amscd}
\usepackage{mathrsfs}
\usepackage{enumerate}

\usepackage[all]{xy} 
\CompileMatrices

\numberwithin{equation}{section}

\newcommand{\C}{{\mathbb C}}
\newcommand{\Q}{\mathbb {Q}}

\newcommand{\Z}{{\mathbb Z}}
\newcommand{\B}{{\mathcal{B}}}

\newcommand{\gl}{{\mathfrak{gl}}}

\newcommand{\seteq}{\mathbin{:=}}

\theoremstyle{plain}
\newtheorem{lemma}{Lemma}[section]
\newtheorem{prop}[lemma]{Proposition}
\newtheorem{theorem}[lemma]{Theorem}
\newcommand{\Prop}{\begin{prop}}
\newcommand{\enprop}{\end{prop}}
\newcommand{\Lemma}{\begin{lemma}}
\newcommand{\enlemma}{\end{lemma}}
\newcommand{\Th}{\begin{theorem}}
\newcommand{\enth}{\end{theorem}}
\newtheorem{corollary}[lemma]{Corollary}
\newcommand{\Cor}{\begin{corollary}}
\newcommand{\encor}{\end{corollary}}
\newtheorem{definition}[lemma]{Definition}

\newcommand{\Def}{\begin{definition}}
\newcommand{\edf}{\end{definition}}

\theoremstyle{definition}
\newtheorem{remark}[lemma]{Remark}

\newtheorem{conjecture}[lemma]{Conjecture}

\newcommand{\g}{{\mathfrak{g}}}

\newcommand{\isoto}[1][]%
{{\mathop{\buildrel{\sim}\over\longrightarrow}\limits_{#1}}}
\renewcommand{\hom}{\operatorname{\it \mathscr{H}\kern-.25em om}}

\newcommand{\eq}{\begin{eqnarray}}
\newcommand{\eneq}{\end{eqnarray}}

\newcommand{\eqn}{\begin{eqnarray*}}
\newcommand{\eneqn}{\end{eqnarray*}}

\newenvironment{tenumerate}{
  \begin{enumerate}
  
  }{\end{enumerate}}

\newcommand{\bnum}{\begin{tenumerate}}
\newcommand{\enum}{\end{tenumerate}}
\newenvironment{anumerate}{
  \begin{enumerate}
  
  }{\end{enumerate}}

\newcommand{\bia}{\begin{anumerate}}
\newcommand{\eia}{\end{anumerate}}
\newcommand{\on}{\operatorname}

\newcommand{\bni}{\begin{tenumerate}}
\newcommand{\eni}{\end{tenumerate}}

\newcommand{\QED}{\end{proof}}

\newcommand{\soplus}{\mathop{\mbox{\normalsize$\bigoplus$}}\limits}

\newcommand{\cl}{\colon}

\newcommand{\id}{\on{id}}
\newcommand{\ba}{\begin{array}}
\newcommand{\ea}{\end{array}}

\newcommand{\bi}{\bni}
\newcommand{\ei}{\eni}
\newcommand{\mono}{\rightarrowtail}

\newcommand{\set}[2]{\left\{#1 \mathbin{;} #2 \right\}}

\newcommand{\hs}{\hspace}

\newcommand{\eqsub}{\begin{subequations}\begin{eqnarray}}
\newcommand{\eneqsub}{\end{eqnarray}\end{subequations}}

\newcommand{\A}{\mathbf{A}}

\renewcommand{\le}{\leqslant}
\renewcommand{\ge}{\geqslant}
\newcommand{\Ind}{{{\on{Ind}}}}

\newcommand{\nc}{\newcommand}
\nc{\la}{\lambda}
\nc{\lam}{\lambda}
\nc{\U}[1][\g]{U_q(#1)}
\nc{\te}{\tilde{e}}
\nc{\tei}{\tilde{e}_i}
\nc{\tf}{\tilde{f}}
\nc{\tfi}{\tilde{f}_i}
\nc{\tU}{\widetilde U_q(\g)}
\nc{\tE}{\widetilde{E}}
\nc{\tF}{\widetilde{F}}

\nc{\BZ}{{\mathbb{Z}}}
\nc{\al}{\alpha}
\nc{\qs}{{q}}
\nc{\lan}{\langle}
\nc{\ran}{\rangle}
\nc{\re}{{\mathrm{re}}}
\nc{\wt}{\operatorname{wt}}
\nc{\Uf}[1][\g]{U^-_q(#1)}
\nc{\Ue}{U^+_q(\g)}
\nc{\eps}{\varepsilon}
\nc{\vphi}{\varphi}
\nc{\sphi}{\varphi^*}
\nc{\seps}{\varepsilon^*}

\nc{\nn}{\nonumber}

\nc{\vp}{\varpi}

\nc{\cls}{{\operatorname{cl}}}

\nc{\Wt}{{\operatorname{Wt}}}
\nc{\Us}{U'_q(\g)}
\nc{\La}{\Lambda}
\nc{\ro}{{\rm(}}
\nc{\rf}{{\rm)}}
\nc{\norm}{{\mathrm{norm}}}
\nc{\qbox}{\quad\mbox}
\nc{\braid}{{\mathfrak{B}}}
\nc{\Ad}{\operatorname{Ad}}
\nc{\Aut}{\operatorname{Aut}}
\nc{\dt}[1]{\tilde{\tilde #1}}
\nc{\Sn}{S^{{\mathrm{norm}}}}
\nc{\aff}{{\mathrm{aff}}}
\nc{\rk}{{\mathrm{rk}}}
\nc{\tQ}{\widetilde{Q}}
\nc{\tP}{\widetilde{P}}
\nc{\tW}{\widetilde{W}}
\nc{\Dyn}{\mathrm{Dyn}}
\nc{\tD}{\widetilde{\Delta}}
\nc{\height}{{\operatorname{ht}}}
\nc{\bl}{\bigl}
\nc{\br}{\bigr}
\nc{\Hecke}{\mathrm{H}}
\nc{\HA}{\Hecke^{\mathrm{A}}}
\nc{\HB}{\Hecke^{\mathrm{B}}}
\nc{\HD}{\Hecke^{\mathrm{D}}}
\nc{\K}{\mathbf{K}}
\newcommand{\scbul}{{\,\raise1pt\hbox{$\scriptscriptstyle\bullet$}\,}}
\nc{\vac}{{\phi}}
\nc{\Bt}{\B_{\theta}(\g)}
\nc{\Vt}{V_{\theta}}
\nc{\Lt}{L_{\theta}}
\nc{\CBt}{B_{\theta}}
\nc{\be}{\begin{enumerate}}
\nc{\ee}{\end{enumerate}}
\nc{\low}{{\mathrm{low}}}
\nc{\upper}{{\mathrm{up}}}
\nc{\Zodd}{\Z_{\mathrm{odd}}}
\nc{\Zeven}{\Z_{\mathrm{even}}}
\nc{\Ft}[1][n]{\mathbb{P}\mathrm{ol}_{#1}}
\nc{\Ftf}[1][n]{\widetilde{\mathbb{P}\mathrm{ol}}_{#1}}
\nc{\KA}{\on{K}^{\mathrm{A}}}
\nc{\KB}{\on{K}^{\mathrm{B}}}
\nc{\KD}{\on{K}^{\mathrm{D}}}
\nc{\Res}{\on{Res}}
\nc{\Fc}[1][{n,m}]{\mathbf{F}_{#1}}
\nc{\tphi}{\tilde{\varphi}}
\nc{\CO}{\mathscr{O}}
\nc{\CS}{\mathscr{S}}

\begin{document}

\title
{Crystals and
affine Hecke algebras of type D}

\author{Masaki KASHIWARA}
\address{Research Institute for Mathematical Sciences,
Kyoto University, Kyoto 606, Japan
}
\thanks{The first author is partially supported by 
Grant-in-Aid for Scientific Research (B) 18340007,
Japan Society for the Promotion of Science.}

\author{Vanessa MIEMIETZ}
\address{Mathematisches Institut, Universit\"at zu K\"oln, 50931 K\"oln, Germany}
%\thanks{The second author is supported by }

%\date{\today}

\begin{abstract}
The Lascoux-Leclerc-Thibon-Ariki theory
asserts that the K-group of the representations of
the affine Hecke algebras of type A is isomorphic to
the algebra of functions
on the maximal unipotent subgroup of the group associated with a Lie algebra
$\g$ where $\g$ is $\gl_\infty$ or 
the affine Lie algebra $A^{(1)}_\ell$,
and the irreducible representations correspond to the
upper global bases.
Recently, N. Enomoto and the first author 
presented the notion of symmetric crystals and
formulated analogous conjectures for
the affine Hecke algebras of type B.
In this note, we present similar conjectures
for certain classes of irreducible representations
of affine Hecke algebras of type D.
The crystal for type D is a double cover of 
the one for type B.
\end{abstract}

\maketitle

\section{Introduction}
Lascoux-Leclerc-Thibon (\cite{LLT})
conjectured the relations between the representations of 
Hecke algebras of {\em type A} and the crystal bases of the affine Lie algebras
of type A. Then, S.~Ariki (\cite{A}) observed that it should be understood 
in the setting of affine Hecke algebras
and proved the LLT conjecture in a more general framework.
Recently, N. Enomoto and the first author 
presented the notion of symmetric crystals and
conjectured that 
certain classes of irreducible representations 
of the affine Hecke algebras of {\em type B} are described by
symmetric crystals (\cite{EK}).

The purpose of this note is to
formulate and explain 
conjectures on certain classes of irreducible representations
of affine Hecke algebras of {\em type D} and symmetric crystals.

Let us begin by recalling the Lascoux-Leclerc-Thibon-Ariki theory.
Let $\HA_n$ be the affine Hecke algebra of type A of degree $n$.
Let $\KA_n$ be the Grothendieck group of
the abelian category of finite-dimensional $\HA_n$-modules,
and $\KA=\soplus\nolimits_{n\ge0}\KA_n$.
Then it has a structure of Hopf algebra by the restriction and
the induction functors. %(cf.\ \S \ref{subsec:affA}).
The set $I=\C^*$ may be regarded as a Dynkin diagram
with $I$ as the set of vertices and with edges between $a\in I$ and $ap^2$.
%(see \eqref{eq:Dy}).
Here $p$ is the parameter of 
the affine Hecke algebra, usually denoted by $q$.
Let $\g_I$ be the associated Lie algebra, and
$\g_I^-$ the unipotent Lie subalgebra. Let $U_I$ be the group associated
to $\g_I^-$.
Hence $\g_I$ is isomorphic to a direct sum of copies of $A^{(1)}_\ell$
if $p^2$ is a primitive $\ell$-th root of unity
and to a direct sum of copies of $\gl_\infty$ if $p$ has an infinite order.
Then $\C\otimes \KA$ is isomorphic to the algebra $\CO(U_I)$
of regular functions on $U_I$.
Let $\U[\g_I]$ be the associated quantized enveloping algebra.
Then $\Uf[\g_I]$ has an upper global basis 
$\{G^{\upper}(b)\}_{b\in B(\infty)}$.
By specializing $\soplus \C[q,q^{-1}]G^\upper(b)$
at $q=1$, we obtain $\CO(U_I)$.
Then the LLTA theory says that the
elements associated to irreducible $\HA$-modules
corresponds to the image of the upper global basis.

In \cite{EK}, N.~Enomoto and the first author
gave analogous conjectures for affine Hecke algebras of type B.
In the type B case, we have to replace 
$\Uf[\g_I]$ and its upper global basis with
a new object, the symmetric crystals.
It is roughly stated as follows.
Let $\HB_n$ be the affine Hecke algebra of type B of degree $n$.
Let $\KB_n$ be the Grothendieck group of
the abelian category of finite-dimensional modules over $\HB_n$,
and $\KB=\oplus_{n\ge0}\KB_n$.
Then $\KB$ has a structure of a Hopf bimodule over $\KA$.
The group $U_I$ has the anti-involution $\theta$ induced by
the involution $a\mapsto a^{-1}$ of $I=\C^*$. Let $U_I^\theta$
be the $\theta$-fixed point set of $U_I$.
Then $\CO(U_I^\theta)$ is a quotient ring of $\CO(U_I)$.
The action of $\CO(U_I)\simeq\C\otimes\KA$ on $\C\otimes\KB$,
in fact, descends to the action of $\CO(U_I^\theta)$.
They introduced
the algebra $\Bt$, a kind of a $q$-analogue 
of the ring of differential operators on $U_I^\theta$
and then $V_\theta(\la)$, a $q$-analogue 
of $\CO(U_I^\theta)$.
The module $V_\theta(\la)$ is an irreducible $\Bt$-module
generated by the highest weight vector $\vac_\la$.
Then they conjectured that:
\bi
\item
$V_\theta(\la)$ has a crystal basis and an upper global basis.
\item
$\KB$ is isomorphic to a specialization of $V_\theta(\la)$ at $q=1$
as an $\CO(U_I)$-module, and
the irreducible representations correspond
to the upper global basis of $V_\theta(\la)$ at $q=1$.
\ei
The representations of $\HB_n$
such that some of $X_i$ have an eigenvalue $\pm1$ are excluded.

In this note, we treat the affine Hecke algebras of type D.
Let
$\HD_n$ be the affine Hecke algebra of type D of degree $n$
($\HD_0=\C\oplus\C$, $\HD_1=\C[X_1^\pm]$, see \S\;\ref{subsec:def}).
Let $\KD_n$ be the Grothendieck group of finite-dimensional $\HD_n$-modules,
and set $\KD=\soplus\nolimits_{n\ge0}\KD_n$.
In D-case, we use the same algebra $\Bt$,
but, instead of $V_\theta(\la)$,
we use a $\Bt$-module $\Vt$ generated
by a pair of highest weight vectors
$\vac_\pm$ (see \S\;\ref{subsec:symmetry}). Our conjecture 
(see \S\;\ref{subsec:main}) is then:
\bi
\item
$\Vt$ has a crystal basis and an upper global basis.
\item
$\KD$ is isomorphic to a specialization of $\Vt$ at $q=1$, 
and the irreducible representations correspond
to the upper global basis of $\Vt$ at $q=1$.
\ei
The representations of $\HD_n$
such that some of $X_i$ have an eigenvalue $\pm1$ are again excluded.

Note that the crystal basis for type D is a double cover of the one for
type B.

\section{Symmetric crystals}\label{sec:crystal}

\subsection{Quantized universal enveloping algebras}
We shall recall the quantized universal enveloping algebra
$U_q(\g)$.
Let $I$ be an index set (for simple roots),
and $Q$ the free $\Z$-module with a basis $\{\al_i\}_{i\in I}$.
Let $(\scbul,\scbul)\cl Q\times Q\to\Z$ be
a symmetric bilinear form such that
$(\al_i,\al_i)/2\in\Z_{>0}$ for any $i$ and
$(\al_i^\vee,\al_j)\in\Z_{\le0}$ for $i\not=j$ where
$\al_i^\vee\seteq2\al_i/(\al_i,\al_i)$.
Let $q$ be an indeterminate and set 
$\K\seteq\Q(\qs)$.
We define its subrings $\A_0$, $\A_\infty$ and $\A$ as follows.
\eq
&&\ba{rcl}
\A_0&=&\set{f/g}{f(\qs),g(\qs)\in \Q[\qs],\, g(0)\not=0},\\[3pt]
\A_\infty&=&\set{f/g}{f(\qs^{-1}),g(\qs^{-1})\in \Q[\qs^{-1}],\, g(0)\not=0},
\\[3pt]
\A&=&\Q[\qs,\qs^{-1}].
\ea
\eneq

\begin{definition}\label{U_q(g)}
The quantized universal enveloping algebra $U_q(\g)$ is the $\K$-algebra
generated by the elements $e_i,f_i$ and invertible
elements $t_i\ (i\in I)$
with the following defining relations.
\begin{enumerate}[{\rm(1)}]

\item The $t_i$'s commute with each other.

\item
$t_je_i\,t_j^{-1}=q^{(\al_j,\al_i)}\,e_i\ $
and $\ t_jf_it_j^{-1}=q^{-(\al_j,\al_i)}f_i\ $
for any $i,j\in I$.

\item\label{even} $\lbrack e_i,f_j\rbrack
=\delta_{ij}\dfrac{t_i-t_i^{-1}}{q_i-q_i^{-1}}$
for $i$, $j\in I$, where $q_i\seteq q^{(\al_i,\al_i)/2}$.

\item {\rm(}{\em Serre relation}{\rm)} For $i\not= j$,
\begin{eqnarray*}
&&\sum^b_{k=0}(-1)^ke^{(k)}_ie_je^{(b-k)}_i=\sum^b_{k=0}(-1)^kf^{(k)}_i
f_jf_i^{(b-k)}=0.
\end{eqnarray*}
Here $b=1-(\al_i^\vee,\al_j)$ and
\begin{eqnarray*}
e^{(k)}_i=e^k_i/\lbrack k\rbrack_i!\ ,&& f^{(k)}_i=f^k_i/\lbrack k\rbrack_i!\ ,\\
\lbrack k\rbrack_i=(q^k_i-q^{-k}_i)/(q_i-q^{-1}_i)\ ,
&&\lbrack k\rbrack_i!=\lbrack 1\rbrack_i\cdots \lbrack k\rbrack_i\,.
\end{eqnarray*}
\end{enumerate}
\end{definition}

Let us denote by $\Uf$ (resp.~$\Ue$)
the subalgebra of $\U$ generated by the $f_i$'s (resp.{} the $e_i$'s).
Let $e'_i$ and $e^*_i$ be the operators on $\Uf$ defined by
$$[e_i,a]=\dfrac{(e^*_ia)t_i-t_i^{-1}e'_ia}{q_i-q_i^{-1}}\quad(a\in\Uf).$$
Then these operators satisfy the following formula similar to derivations:
$$e_i'(ab)=e_i'(a)b+(\Ad(t_i)a)e_i'b,\quad
e_i^*(ab)=ae_i^*b+(e_i^*a)(\Ad(t_i)b).$$
The algebra $\Uf$ has a unique symmetric bilinear form $(\scbul,\scbul)$
such that $(1,1)=1$ and
$$(e'_ia,b)=(a,f_ib)\quad\text{for any $a,b\in\Uf$.}$$
It is non-degenerate and satisfies $(e^*_ia,b)=(a,bf_i)$.

\subsection{Symmetry}\label{subsec:symmetry}
Let $\theta$ be an automorphism of
$I$ such that $\theta^2=\id$ and 
$(\al_{\theta(i)},\al_{\theta(j)})=(\al_i,\al_j)$.
Hence it extends to an automorphism of the root lattice $Q$
by $\theta(\al_i)=\al_{\theta(i)}$,
and induces an automorphism of $\U$.

Let $\Bt$ be the $\K$-algebra
generated by $E_i$, $F_i$, and
invertible elements $K_i$ ($i\in I$)
satisfying the following defining relations:
\eq&&\left\{
\parbox{30em}{
\bi
\item the $K_i$'s commute with each other,
\item
$K_{\theta(i)}=K_i$ for any $i\in I$,
%\item
%$\sigma^2=1$
\item
$K_iE_jK_i^{-1}=q^{(\al_i+\al_{\theta(i)},\al_j)}E_j$ and
$K_iF_jK_i^{-1}=q^{(\al_i+\al_{\theta(i)},-\al_j)}F_j$
for $i,j\in I$,
\item
$E_iF_j=q^{-(\al_i,\al_j)}F_jE_i+
(\delta_{i,j}+\delta_{\theta(i),j}K_i)$
\quad for $i,j\in I$,
\item
the $E_i$'s and the $F_i$'s satisfy the Serre relations.
%\item $\sigma$ commutes with $E_i,K_i,F_i$ for $i \in I$
\ei}
\right.
\label{rel:EFK}
\eneq
Hence $\Bt\simeq\Uf\otimes\K[K_i^{\pm1};i\in I]\otimes \Ue$.
We set $E_i^{(n)}=E_i^n/[n]_i!$
and $F_i^{(n)}=F_i^n/[n]_i!$. 
%We define $\Uf_\A$ as the $\A$-subalgebra of $\Bt$ generated 
%by $F_i^{(n)}$ ($i \in I$).

%Let $\la\in P_+\seteq\set{\la\in \Hom(Q,\Q)}%
%{\mbox{$\lan \al_i^\vee,\la\ran\in\Z_{\ge0}$ for any $i\in I$}}$ 
%be a dominant integral weight such that
%$\theta(\la)=\la$.
\Prop\label{prop:Vtheta}
\bi
\item
There exists a $\Bt$-module $\Vt$
generated by  linearly independent vectors $\vac_+$ and $\vac_-$ such that
\be[{\rm(a)}]
\item
$E_i\vac_\pm=0$ for any $i\in I$,
\item
$K_i\vac_\pm=\vac_\mp$ for any $i\in I$,
\item
$\set{u\in \Vt}{\text{$E_iu=0$ for any $i\in I$}}
=\K\vac_+ \oplus \K\vac_-$.
\ee
Moreover such a $\Vt$ is unique up to an isomorphism.
\item
There exists a unique symmetric bilinear form $(\scbul,\scbul)$
on $\Vt$ such that 
$(\vac_{\eps_1},\vac_{\eps_2})
=\delta_{\eps_1,\eps_2}$ for $\eps_1,\eps_2 \in \{+,-\}$ and
$(E_iu,v)=(u,F_iv)$ for any $i\in I$ and $u,v\in \Vt$,
and it is non-degenerate.
\ei
\enprop

%The pair $(\Bt,\Vt)$ is an analogue of $(\B,U_q^-(\g))$.
Such a $\Vt$ is constructed as follows.

Let $\CS$ be the quantum shuffle algebra 
(see \cite{Leclerc}) generated by words 
$\lan i_1, \dots,i_l \ran$ for $i_1, \dots,i_l \in I$ and $l \geq 1$ 
and $\vac_+''$ and $\vac_-''$ as two empty words. 
%This has a left and a right action of $U_q^-(\g)$ 
%by the embedding $f_i \mapsto \lan i \ran$. 
We assign to a word $\lan i_1, \dots, i_l \ran $ 
the weight $-(\al_{i_1}+ \cdots +\al_{i_l})$.
We define the actions of $E_i$, $F_i$
and $K_i$ on $\CS$ as follows:
\eq&&\ba{rl}
F_i \vac''_+&= \lan i \ran,\quad F_i \vac''_-= \lan \theta i \ran,\\[2pt]
E_i \lan j \ran &= \delta_{i,j}\vac''_+ +\delta_{i,\theta j}\vac_-'',\\[2pt]
K_i \vac''_\pm &= \vac''_\mp,\\[2pt]
K_i\lan i_1, \dots, i_l \ran&=
q^{-(\al_{i}+\al_{\theta(i)},\al_{i_1}+\cdots+\al_{i_{l}})}
\lan i_1,\dots,i_{l-1}, \theta(i_l) \ran,\\[2pt]
%K_i a \vac''_\pm &= (\Ad(t_it_{\theta(i)})a)\vac''_\pm\\[2pt]
E_i\lan i_1, \dots, i_l \ran&= \delta_{i,i_1}\lan i_2,\dots, i_l \ran,\\[3pt]
%F_i(\lan i_1, \dots, i_l \ran )&=
%f_i\lan i_1, \dots, i_l \ran + q^{(\al_i, -(\al_{i_1}+ \cdots +\al_{i_{l-1}} + \al_{\theta i_l} ))}  \lan i_1, \dots, i_{l-1}, \theta i_{l} \ran  f_{\theta(i)}.
F_i \lan i_1,\ldots,i_l\ran&=
\lan i \ran *\lan i_1,\ldots,i_l\ran+q
^{(\al_i,\wt(\lan i_1,\ldots,i_{l-1},\theta(i_l)
\ran))}\lan i_1,\ldots,i_{l-1},\theta(i_l)\ran * \lan \theta i \ran\\[2pt]
&=\sum\limits_{\nu=0}^lq^{-(\al_i,\al_{i_1}+\cdots+\al_{i_\nu})}
\lan i_1,\ldots,i_\nu,i,i_{\nu+1},\ldots,i_l\ran\\
&\quad+ q^{-(\al_i,\al_1+\cdots +\al_{i_{l-1}}+\al_{\theta(i_l)})}
\sum\limits_{\nu=0}^l
q^{-(\al_{\theta(i)},\al_{\nu+1}+\cdots+\al_{i_{l-1}}+\al_{\theta(i_l)})}\\
&\hs{12em}
\cdot\lan i_1,\ldots,i_\nu,\theta(i),i_{\nu+1},\ldots,i_{l-1},\theta(i_l)\ran
\ea
\eneq
for $i,j \in I$, $l\ge1$ and $i_1,\ldots,i_l\in I$.

Then the operators $E_i$, $F_i$ and $K_i$ satisfy the commutation relations
\eqref{rel:EFK} except the Serre relations for the $E_i$'s.

Consider the $U_q^-(\g)$-module 
$V'=U_q^-(\g) \vac'_+ \oplus U_q^-(\g) \vac'_-$ generated by a pair of vacuum vectors $\vac'_\pm$. 
There exists a unique $U_q^-(\g)$-linear map 
$\psi\cl V' \rightarrow \CS$ such that $\vac'_\pm \mapsto \vac''_\pm$. 
We define an action of $\Bt$ on $V'$ by 
\eq&&\ba{rcl}
%\sigma \vac'_+ &=& \vac'_-\\[2pt]
K_i (a \vac'_\pm) &=& (\Ad(t_it_{\theta(i)})a)\vac'_\mp\\[2pt]
E_i (a \vac'_\pm) &=& e_i'(a)\vac'_\pm + Ad(t_i)(e_{\theta i}^*(a)) \vac'_\mp\\[2pt]
F_i (a \vac'_\pm)&=& f_ia \vac'_\pm.
\ea
\qquad\text{for $a\in\Uf$.}
\eneq
Then $\psi$ commutes with the actions of $E_i$, $F_i$ and $K_i$,
and
its image $\psi(V')$ is $\Vt$.

Hereafter we assume further that
\eq&&
\text{there is no $i\in I$ such that $\theta(i)=i$.}
\eneq
Under this condition, we conjecture that $\Vt$ has a crystal basis.
This means the following.
We define the modified root operators:
$$\tE_i(u)=\sum_{n\ge1}F_i^{(n-1)}u_n\quad
\text{and}\quad
\tF_i(u)=\sum_{n\ge0}F_i^{(n+1)}u_n
$$
when writing $u=\sum_{n\ge0}F_i^{(n)}u_n$ with $E_iu_n=0$.
Let $\Lt$ be the $\A_0$-submodule of
$\Vt$ generated by $\tF_{i_1}\cdots\tF_{i_\ell}\vac_\pm$
($\ell\ge0$ and $i_1,\ldots,i_\ell\in I$\,),
and define the subset $\CBt\subset \Lt/\qs \Lt$ 
by:
$$\CBt\seteq\set{\tF_{i_1}\cdots\tF_{i_\ell}\vac_\pm\bmod \qs \Lt}%
{\text{$\ell\ge0$, $i_1,\ldots, i_\ell\in I$}}.$$

\begin{conjecture}\label{conj:crystal}
\bi
\item
$\tF_i\Lt\subset \Lt$
and $\tE_i\Lt\subset \Lt$,
\item
$\CBt$ is a basis of $\Lt/\qs\Lt$,
\item
$\tF_i\CBt\subset \CBt$,
and
$\tE_i\CBt\subset \CBt\sqcup\{0\}$.
\ei
\end{conjecture}

Moreover we conjecture that
$\Vt$ has a global crystal basis.
Namely, let $-$ be the bar-operator of
$\Vt$, which is characterized by:
$\overline{q}=q^{-1}$, $-$ commutes with the $E_i$'s, and
$(\vac_\pm)^-=\vac_\pm$
% commute given by $-\cl a\vac_\pm\to \bar a\vac_?$ ($a\in \Uf$)
(such an operator exists).
Let us denote by $\Bt_\A^\upper$ the $\A$-subalgebra of $\Bt$ generated 
by $E_i^{(n)}$, $F_i$ and $K_i^{\pm1}$ ($i \in I$).
Let $(\Vt)_\A$ be the largest $\Bt^\upper_\A$-submodule
of $\Vt$ such that $(\Vt)_\A\cap(\K\vac_++\K\vac_-)=\A\vac_++\A\vac_-$.

\begin{conjecture}\label{conj:balanced}
$(\Lt,\Lt^-,(\Vt)_\A)$
is balanced.
\end{conjecture}
Namely,
$E\seteq \Lt\cap\Lt^-\cap(\Vt)_\A\to
\Lt/\qs\Lt$ is an isomorphism.
Let $G^\upper\colon \Lt/\qs \Lt\isoto 
E$ be its inverse.
Then $\set{G^\upper(b)}{b\in \CBt}$ forms a basis of $\Vt$.
We call this basis the {\em upper global basis} of $\Vt$.

\begin{remark}\label{rem:crystal}
Assume that Conjectures \ref{conj:crystal} and \ref{conj:crystal}
hold.
\be[(i)]
\item We have
$\set{b\in\CBt}{\text{$\tE_ib=0$ for any $i\in I$}}
=\{\vac_+,\vac_-\}$.
\item
There exists a unique involution
$\sigma$ of the $\Bt$-module $\Vt$ such that
$\sigma(\vac_\pm)=\vac_\mp$.
It extends to the involution $\sigma$ of $\CS$ by
$\sigma(\lan i_1,\ldots,i_l\ran)=\lan i_1,\ldots,i_{l-1},\theta(i_l)\ran$.
It induces also involutions of $\Lt$ and $\CBt$.
\item
We have $\sigma(b)\not=b$ for any $b\in\CBt$.
\item
We conjecture that $\tF_ib\not=\tF_jb$ for
any $b\in\CBt$ and $i\not=j\in I$.
\item
In \cite{EK}, a $\Bt$-module $V_\theta(\la)=\Bt\vac_\la$
and its crystal basis $B_\theta(\la)$ are introduced. 
We have a monomorphism of $\Bt$-modules
$$V_\theta(\la)\mono \Vt$$
with $\la=0$, which
sends $\vac_\la$ to $\vac_++\vac_-$. Its image
coincides with $\set{v\in\Vt}{\sigma(v)=v}$
Any element $b\in B_\theta(\la)$
is sent to $b'+\sigma(b')$ for some $b'\in\CBt$.
Moreover, we have $\iota(G^\upper(b))=G^\upper(b')+\sigma (G^\upper(b'))$.
In particular, we have
$$B_\theta(\la)\simeq\CBt/\sim.$$
Here $\sim$ is the equivalence relation given by $b\sim\sigma b$.
\label{rem:BD}
\ee
\end{remark}

\section{Affine Hecke algebra of type D}\label{sec:aff}
\subsection{Definition}\label{subsec:def}

For $p \in\C^*$ and $n\in\Z_{\ge2}$, 
the affine Hecke algebra $\HD_n$ of type $D_n$
is the $\C$-algebra generated by
$T_i$ ($0\le i<n$) and invertible elements $X_i$ ($1\le i\le n$)
satisfying the defining relations:
\bi
\item
the $X_i$'s commute with each other,
\item
the $T_i$'s satisfy the braid relation:
$T_1T_0=T_0T_1$,
$T_0T_2T_0=T_2T_0T_2$,
$T_iT_{i+1}T_i=T_{i+1}T_iT_{i+1}$ ($1\le i<n-1$),
$T_iT_j=T_jT_i$ ($1\le i<j-1<n-1$ or $i=0<3\le j<n$),
\item
$(T_i-p)(T_i+p^{-1})=0$ ($0\le i<n$),
\item
$T_0X_1^{-1}T_0=X_2$,
$T_iX_iT_i=X_{i+1}$ ($1\le i<n$),
and $T_iX_j=X_jT_i$ if $1 \leq i\not=j,j-1$ or $i=0$ and $j \geq 3$.
\ei
We define $\HD_0=\C \oplus \C$ and $\HD_1=\C[X_1^{\pm1}]$.

We assume that $p\in\C^*$ satisfies
\eq
p^2\not=1.
\eneq
Let us denote by $\Ft$ the Laurent polynomial ring
$\C[X_1^{\pm1},\ldots,X_n^{\pm1}]$,
and by $\Ftf$ its quotient field $\C(X_1,\ldots, X_n)$.
Then $\HD_n$ is isomorphic to the tensor product
of $\Ft$ and the subalgebra generated by the $T_i$'s
that is isomorphic to the Hecke algebra of type $D_n$.
We have
$$
T_ia=(s_ia)T_i+(p-p^{-1})\dfrac{a-s_ia}{1-X^{-\al_i^\vee}}
\quad\text{for $a\in\Ft$.}$$
Here, $X^{-\al_i^\vee}=X_1^{-1}X_2^{-1}$ ($i=0$)
and $X^{-\al_i^\vee}=X_{i}X_{i+1}^{-1}$ ($1\le i<n$).
The $s_i$'s are the Weyl group action on $\Ft$:
$(s_0a)(X_1,\ldots,X_n)=a(X_2^{-1},X_1^{-1},\ldots,X_n)$ 
and
$(s_ia)(X_1,\ldots,X_n)=a(X_1,\ldots,X_{i+1},X_i,\ldots,X_n)$ for $1\le i<n$.

\subsection{Intertwiner}
The algebra $\HD_n$ acts faithfully on
$\HD_n/\sum_i\HD_n(T_i-p)\simeq\Ft$.
Set $\vphi_i=(1-X^{-\al_i^\vee})T_i-(p-p^{-1})\in\HD_n$
and $\tphi_i=(p^{-1}-pX^{-\al_i^\vee})^{-1}\vphi_i\in\Ftf\otimes_{\Ft}
\HD_n$.
Then the action of $\tphi_i$ on $\Ft$ coincides with $s_i$.
They are called intertwiners.

\subsection{Affine Hecke algebra of type A}\label{subsec:affA}
The affine Hecke algebra $\HA_n$ of type $A_n$
is isomorphic to the subalgebra
of $\HD_n$ generated by $T_i$ ($1\le i<n$) and $X_i^{\pm1}$
($1\le i\le n$). For a finite-dimensional $\HA_n$-module $M$,
let us decompose
\eq\label{eq:wtdec}
M=\soplus_{a\in(\C^*)^n}M_a
\eneq
where $M_a=\set{u\in M}{\text{$(X_i-a_i)^Nu=0$ for any $i$ %=1,\ldots,n$
and $N\gg0$}}$
for $a=(a_1,\ldots,a_n)\in(\C^*)^n$.
For a subset $I\subset\C^*$,
we say that $M$ is of type $I$ if all the eigenvalues of $X_i$
belong to $I$.
The group $\Z$ acts on $\C^*$ by $\Z\ni n\cl a\mapsto ap^{2n}$.
By well-known results in type A,
it is enough to treat the irreducible modules of type $I$
for an orbit $I$ with respect to the $\Z$-action on $\C^*$ in order to study the irreducible modules over the affine Hecke
algebras of type A.

\subsection{Representations of affine Hecke algebras of type D}
\label{subsec:main}

For $n,m\ge 0$, set
$\Fc\seteq\C[X_1^{\pm1},\ldots,X_{n+m}^{\pm1},
D^{-1}]$
where
$$D\seteq\prod\limits_{1\le i\le n<j\le n+m}\hs{-2ex}
(X_i-p^2 X_j)(X_i-p^{-2}X_j)
(X_i-p^2 X_j^{-1})(X_i-p^{-2}X_j^{-1})(X_i-X_j)(X_i-X_j^{-1}).$$
Then we can embed
$\HD_m$ into $\HD_{n+m}\otimes_{\Ft[n+m]}\Fc$
by 
\eqn
T_0&\mapsto& \tphi_{n}\cdots\tphi_1\tphi_{n+1}\cdots\tphi_2T_0\tphi_2\cdots\tphi_{n+1}\tphi_1\cdots\tphi_n,\\ 
T_i&\mapsto& T_{i+n}\quad(1\le i<m),\\ 
X_i&\mapsto& X_{i+n}\quad(1\le i\le m).
\eneqn
Its image commutes with $\HD_n\subset\HD_{n+m}$.
Hence $\HD_{n+m}\otimes_{\Ft[n+m]}\Fc$
is a right $\HD_n\otimes\HD_m$-module.

For a finite-dimensional $\HD_n$-module $M$,
we decompose $M$ as in \eqref{eq:wtdec}.
The semidirect product group
$\Z_2\times\Z=\{1,-1\}\times\Z$ acts on $\C^*$
by $(\epsilon,n)\cl a\mapsto a^\epsilon p^{2n}$.

Let $I$ and $J$ be $\Z_2\times\Z$-invariant subsets of
$\C^*$ such that $I\cap J=\emptyset$.
Then for an $\HD_n$-module $N$
of type $I$ and $\HD_m$-module $M$ of type $J$,
the action of $\Ft[n+m]$ on $N\otimes M$ extends to an action of $\Fc$.
We set 
$$N\diamond M\seteq
(\HD_{n+m}\otimes_{\Ft[n+m]}\Fc)
\mathop\otimes\limits
_{(\HD_n\otimes\HD_m)\otimes_{\Ft[n+m]}{\Fc}}(N\otimes M).$$

\Lemma
\bi
\item
Let $N$ be an irreducible $\HD_n$-module
of type $I$ and
$M$ an irreducible $\HD_m$-module of type $J$.
Then $N\diamond M$ is an irreducible $\HD_{n+m}$-module
of type $I\cup J$.
\item Conversely if $L$
is an irreducible $\HD_{n}$-module
of type $I\cup J$, then there exists an integer $m$ $(0\le m\le n)$,
an irreducible $\HD_m$-module $N$
of type $I$ and
an irreducible $\HD_{n-m}$-module $M$ of type $J$
such that $L\simeq N\diamond M$.
\item
Assume that a $\Z_2\times\Z$-orbit $I$ decomposes into
$I=I_+\sqcup I_-$ where $I_\pm$ are $\Z$-orbits and $I_-=(I_+)^{-1}$.
%Assume that $\pm1\not\in I$.
Then for any irreducible $\HD_n$-module $L$ of type $I$,
there exists an irreducible $\HA_n$-module $M$ such that
$L\simeq\Ind_{\HA_n}^{\HD_n}M$.
\ei
\enlemma
Hence in order to study $\HD$-modules, it is enough to
study irreducible modules of type $I$ for a $\Z_2\times\Z$-orbit
$I$ in $\C^*$ such that $I$ is a $\Z$-orbit,
namely $I=\pm\set{p^n}{n\in\Zodd}$ or $I=\pm\set{p^n}{n\in\Zeven}$.
%contains one of $\pm1$.

For a  $\Z_2\times\Z$-invariant subset $I$ of $\C^*$,
we define $\KD_{I,n}$ to be the Grothendieck group of
the abelian category of
finite-dimensional $\HD_n$-modules of type $I$. 
We set $\KD_I=\soplus\nolimits_{n\ge0}\KD_{I,n}$.
%$\KA_{I,n}$ and $\KA_I$ are defined analogously for affine Hecke algebras of type A.
%Then $\KD_I$ is a (right) Hopf $\KA_I$-bimodule
%by the multiplication and the comultiplication
%$$\mu\cl\KD_{I,n}\times \KA_{I,m}\to \KD_{I,n+m}\quad\text{and}\quad
%\Delta\cl\KD_{I,n}\to\soplus_{i+j=n}\KD_{I,i}\otimes\KA_{I,j}$$
%given by $L\otimes M\mapsto \Ind_{\HD_n\otimes\HA_m}^{\HD_{n+m}}(L\otimes M)$
%and $L\mapsto\Res^{\HD_n}_{\HD_i\otimes\HA_j}L$.
%The module $\KD_0$ has two generators $\vac_\pm$.
%We understand
%\eqn&&\ba{l}
%\mu(\vac_+\otimes M)=\Ind_{\HA_n}^{\HD_n}M\quad\parbox[t]{30em}%
%{where $\HA_n$ is identified with the subalgebra\\
%\hs*{10em}$\lan T_1,\ldots,T_{n-1},X_1^\pm,\ldots,X_n^\pm\ran$,}\\[3ex]
%\mu(\vac_-\otimes M)=\Ind_{\HA_n}^{\HD_n}M
%\quad\parbox[t]{30em}%
%{where $\HA_n$ is identified with the subalgebra\\
%\hs*{10em}$\lan T_0,T_2,\ldots,T_{n-1},X_1^\mp,X_2^\pm,\ldots,X_n^\pm\ran$.}
%\ea
%\eneqn

We take the case
$$I=\set{p^n}{n\in\Zodd}$$ and assume that any of $\pm1$
is not contained in $I$.
The set $I$ may be regarded as the set of vertices of a Dynkin diagram.
Let us define an automorphism $\theta$ of $I$ by $a\mapsto a^{-1}$.
Let $\g_I$ be the associated Lie algebra
($\g_I$ is isomorphic to $\gl_\infty$ if $p$ has an infinite order,
and 
isomorphic to $A^{(1)}_\ell$
if $p^2$ is a primitive $\ell$-th root of unity).

For a finite-dimensional $\HD_n$-module $M$ and $a\in I$,
let $E_aM$ be the generalized {$a$-eigenspace} of $X_n$ on $M$,
regarded as an $\HD_{n-1}$-module. 
Let $F_aM$ be the $\HD_{n+1}$-module
$\Ind_{\HD_n\otimes\C[X_{n+1}^\pm]}^{\HD_{n+1}}
(M \otimes (a))$ 
where $(a)$ is the $1$-dimensional representation of 
$\C[X_{n+1}^{\pm 1}]$ on which $X_{n+1}$ acts as $a$.
Then $E_a$ and $F_a$ are exact functors and define
$E_a\cl \KD_{I,n}\to\KD_{I,n-1}$ and
$F_a\cl \KD_{I,n}\to\KD_{I,n+1}$.

For an irreducible $M \in \KD_{I,n}$ and $a\in I$,
define $\tilde e_{a} M \in \KD_{I,n-1}$ 
to be the socle of $E_aM$.
Define $\tilde f_a M \in \KD_{I,n+1}$ 
to be the cosocle of $F_aM$.
In fact, $\tilde f_aM$ is always irreducible, and
$\tilde e_aM$ is a zero module or irreducible.

The ring $\HD_0=\C\oplus\C$ has two irreducible modules $\vac_\pm$.
We understand $$E_a((b))=
\te_a((b))=\begin{cases}\vac_\pm&\text{if $a=b^{\pm1}$,}\\
0&\text{otherwise,}\end{cases}
\qquad\text{and $F_a(\vac_\pm)=\tf_a(\vac_\pm)=(a^{\pm1})$.}$$

Let $\Vt$ be as in Proposition~\ref{prop:Vtheta}.
\begin{conjecture}
\be[(i)]
\item $\KD$ is isomorphic to  $(\Vt)_\A/(q-1)(\Vt)_\A$.
\item
$\Vt$ has a crystal basis and an upper global basis.
\label{(i)}
\item
The elements of $\KD_I$ associated to irreducible representations
correspond to
the upper global basis of $\Vt$ at $q=1$.
\item The operators $\tilde F_i$ and $\tilde E_i$ 
correspond to $\tilde f_i$ and $\tilde e_i$, respectively.
\ee
\end{conjecture}

Consider $\tilde H = \HD_n \otimes \C[\theta]/(\theta^2-1)$ 
with multiplication $\theta T_1=T_0\theta$, 
$X_1\theta = \theta X_1^{-1}$ 
and $\theta$ commuting with all other generators.
Then $\tilde H$ is isomorphic to 
the specialization of the affine Hecke algebra of type B 
in which the generator for the node corresponding 
to the short root has eigenvalues $\pm 1$. 
This explains why the crystal graph in the above case 
is a double covering of the crystal graph 
for the same $\Z_2\times\Z$-orbit in type $B$.
(See Remark \ref{rem:crystal}\;\eqref{rem:BD}.)

\end{document}